\documentclass[12pt]{article}
\usepackage{theorem, enumerate}

\usepackage{amsmath,amssymb}
\usepackage{graphicx}
\usepackage[colorlinks=true,citecolor=black,linkcolor=black,urlcolor=blue]{hyperref}

\theoremstyle{change}
{\theorembodyfont{\slshape}
\newtheorem{theorem}{Theorem.}[section]
\newtheorem{lemma}[theorem]{Lemma.}
\newtheorem{corollary}[theorem]{Corollary.}
\newtheorem{proposition}[theorem]{Proposition.}
\theorembodyfont{\rmfamily}

\def\proof{\noindent{\textsl{Proof. }}}
\def\sqr#1#2{{\vbox{\hrule height.#2pt
    \hbox{\vrule width.#2pt height#1pt \kern#1pt
        \vrule width.#2pt}\hrule height.#2pt}}}
\def\eqed{\sqr53}
\def\qed{%
    \ifmmode\eqno\eqed
    \else\nobreak\ \hfill\eqed\medbreak\fi}

\def\inv{^{-1}}

\def\reals{{\mathbb R}}
\def\ints{{\mathbb Z}}

\def\al{\alpha}
\def\be{\beta}
\def\a{\alpha}
\def\b{\beta}
\def\g{\gamma}
\def\ba{\bar{\alpha}}               
\def\bb{\bar{\beta}}                
\def\bg{\bar{\gamma}}   	

\def\th{\theta}

\def\spn{\mathrm{span}}
             
\def\C#1#2#3{\Gamma_{#1#2#3}}       
\def\N#1{\Gamma_{#1}}               
\def\y#1#2{Y_{#1,#2}}
\def\nomw{{\mathcal N}_{W}}

\newcommand{\Dita}{{Di{\c{t}}{\u{a}}}}
\title{Complex Hadamard Matrices and Strongly Regular Graphs}

\author{Ada Chan\\
\small Department of Mathematics and Statistics\\[-0.8ex]
\small York University\\[-0.8ex] 
\small Toronto, Canada\\
\small\tt ssachan@yorku.ca}

\begin{document}

\maketitle


\begin{abstract}
  Over forty years ago, Goethals and Seidel showed that if the adjacency algebra of
a strongly regular graph $X$ contains a Hadamard matrix then $X$ is either of
Latin square type or of negative Latin square type.
We extend their result to complex Hadamard matrices
and find only three additional families of parameters for which the strongly
regular graphs have complex Hadamard matrices in their adjacency algebras.

Moreover we show that there are only three distance
regular covers of the complete graph that give complex Hadamard matrices.
\end{abstract}

\section{Introduction}

An $n\times n$ complex matrix $W$ is {\sl type~II} if
\begin{equation}
\sum_{x=1}^{n} \frac{W(a,x)}{W(b,x)} = 
\begin{cases}
n & \text{if $a=b$},\\
0 & \text{otherwise.}
\end{cases}
\label{TypeII}
\end{equation}
We say two type~II matrices $W$ and $W'$ are {\sl type~II equivalent} if 
$W'=MWN$ for some invertible monomial matrices $M$ and $N$.

Type II matrices were introduced by Sylvester as inverse orthogonal matrices
in \cite{Sylvester}.
Over a century later, Jones defined a special class of type~II matrices
called spin models to construct link invariants
\cite{MR89m:57005}.
Spin models are matrices that
satisfy three types of conditions corresponding to the three 
Reidemeister moves on link diagrams.
In particular, the Potts model is the $n\times n$ matrix
\begin{equation*}
-u^3I+u\inv(J-I),
\end{equation*}
where $(u^2+u^{-2})^2=n$ and $J$ is the matrix of all ones.
The Potts models give evaluations of Jones polynomial.
Condition~(\ref{TypeII}) corresponds to the second Reidemeister move,
hence the term type~II matrices.

If every entry of an $n\times n$ type~II matrix $W$ has absolute value $1$,
then (\ref{TypeII}) is equivalent to
\begin{equation*}
W \overline{W}^T = nI.
\end{equation*}
In this case, $W$ is a {\sl complex Hadamard matrix}.
Two complex Hadamard matrices are {\sl equivalent} if they are type~II
equivalent.

Complex Hadamard matrices are natural generalization of (real) Hadamard matrices.
Goethals and Seidel determined that the strongly regular graphs,
which contain (real) Hadamard matrices in their adjacency algebras, 
are either of Latin square type and of negative Latin square type \cite{MR0282872}.
In this paper, we show that there are only three more families of parameters
for which the strongly regular graphs give complex Hadamard matrices,
see Theorem~\ref{Thm_FlatSRG}.
In Section~\ref{cover}, we see that the situation is different 
for distance regular covers of the complete graph.  There are only finitely
many such graphs that give complex Hadamard matrices.

Apart from being a generalization of (real) Hadamard matrices,
complex Hadamard matrices have applications in quantum information theory,
operator theory and combinatorics, 
see \cite{MR2244963}, \cite{MR98k:46087} and \cite{MR2460230}. 
In particular, the transition operator of the continuous-time quantum walk of a graph $X$ on $v$ vertices is
$e^{-i t A}$,
where $A$ is the adjacency matrix of $X$.  We say $X$ admits {\sl instantaneous uniform
mixing at time $\tau$} if
\begin{equation*}
\sqrt{v} e^{-i\tau A}
\end{equation*}
is a complex Hadamard matrix.   Using Theorem~\ref{Thm_FlatSRG} in this paper,  Godsil et al.~determined all the strongly regular graphs that admit uniform mixing in \cite{MR3691539}.


Researchers are interested in finding infinite families of complex Hadamard matrices, 
as well as the classification of complex Hadamard matrices of small sizes.
For instance, Ikuta and Munemasa discovered complex Hadamard matrices from a family of 3-class association scheme  \cite{MR3349366}, containing the 
line graph of the Petersen graph which is also mentioned in Section~\ref{cover}.

Given a type-II matrix $W$, Nomura's construction gives the Bose-Mesner algebra $\nomw$ of an association scheme.
In \cite{MR1958007},  Hosoya and Suzuki showed that $\dim \nomw \geq 3$ when $W$ is type II equivalent to a generalized tensor product.
Many known parametric families of 
complex Hadamard matrices \cite{MR2244963}
result from \Dita's construction given in \cite{MR2065675}, which  is a special case of generalized tensor product.
In Section~\ref{Nomura},   we show that $\dim \nomw=2$ for
the complex Hadamard matrices obtained in this paper, with the
exception of a few.  Hence they do not arise from
\Dita's construction.


\section{Strongly Regular Graphs}
\subsection{Preliminary}
We introduce strongly regular graphs, 
please see \cite{MR1002568} and \cite{MR1829620} for details.

A {\sl strongly regular graph} with parameters $(v,k,a,c)$ is 
a $k$-regular graph $X$ with $v$ vertices, that is not complete,
in which every pair of adjacent vertices have exactly $a$ common neighbours
and every pair of distinct non-adjacent vertices have exactly $c$ common neighbours.
The complement of $X$, denoted by $\overline{X}$,
is also strongly regular and it has parameters 
\begin{equation*}
(v, v-k-1, v-2-2k+c,v-2k+a).
\end{equation*}
Moreover $X$ is not connected if and only if it is isomorphic to $mK_{k+1}$
for some $m\geq 2$.  In this case, $a=k-1$ and $c=0$.  
We say a strongly regular graph $X$ is {\sl primitive} if both $X$ and $\overline{X}$ are connected.

Let $A(X)$ denote the adjacency matrix of a strongly regular graph $X$
with parameters $(v,k,a,c)$.  Then $A(X)$ satisfies
\begin{equation*}
A(X)^2-(a-c)A(X)-(k-c)I=cJ.
\end{equation*}
When $a-c=0$, $A(X)$ is the incidence matrix of a symmetric $(v,k,c)$-design.
Similarly when $a-c=-2$,
$A(\overline{X})$ is the incidence matrix of a symmetric
$(v,v-k-1,v-2k+a)$-design.

As $X$ is $k$-regular, $k$ is an eigenvalue of $A(X)$ with eigenvector ${\bf 1}$.
The eigenvalues of $A(X)$ associated with the eigenvectors orthogonal to ${\bf 1}$ satisfy
\begin{equation*}
z^2-(a-c)z-(k-c)=0.
\end{equation*}
Let $\th$ and $\tau$ be the roots of this quadratic equation, where $\th\geq \tau$.
Then $\th+\tau=a-c$ and $\th\tau=-(k-c)\leq 0$ are integers.
It is known that either $X$ has parameters
$(v, \frac{v-1}{2}, \frac{v-5}{4},\frac{v-1}{4})$ 
or both $\th$ and $\tau$ are integers.
In the former case, $X$ is called a {\sl conference graph}.

Note that $A(\overline{X})$ has eigenvalues $v-k-1$, $-\tau-1$ and $-\th-1$.
It follows that if $X$ is primitive then
\begin{equation*}
k>\th >0>\tau 
\quad \text{and}\quad
v=\frac{(k-\th)(k-\tau)}{k+\th\tau}.
\end{equation*}
The multiplicity of $\th$ and $\tau$ are
\begin{equation*}
m_\th = \frac{(v-1)\tau+k}{\tau-\th}
\quad\text{and}\quad
m_\tau = \frac{(v-1)\th+k}{\th-\tau},
\end{equation*}
respectively.
The integrality of $m_\th$ and $m_\tau$ imposes a necessary condition on $(v,k,a,c)$
for the existence of the strongly regular graph with these parameters.

\subsection{A Necessary Condition}
We give all type~II matrices in the adjacency algebra of a strongly regular graph
\cite{TypeIIComb}.
\begin{theorem} 
\label{Thm_SRG}
Let $X$ be a strongly regular graph with $v > 4$ vertices, valency $k$,
and eigenvalues $k$, $\th$ and $\tau$,
where $\th\geq 0 >\tau$.  
Suppose
\[
W := I+xA(X)+yA(\overline{X}).
\]
Then $W$ is a type-II matrix if and only if one of the following holds
\begin{enumerate}[a.]
\item
\label{SRGa}
$y=x =\frac{1}{2}(2-v\pm\sqrt{v^2-4v})$
and $W$ is the Potts model.
\item
\label{SRGb}
$x =\frac{1}{2}(2-v\pm\sqrt{v^2-4v})$ and 
$y = (1+xk)/(x+k)$,
and $X$ is isomorphic to $mK_{k+1}$ for some $m>1$.
\item
\label{SRGc}
$x=1$ and 
\[
y=\frac{v+2(1+\th)(1+\tau) \pm \sqrt{v^2+4v(1+\th)(1+\tau)}}{2(1+\th)(1+\tau)},
\]
and $A(\overline{X})$ is the incidence matrix of a symmetric $(v, v-k-1, \lambda)$-design
where $\lambda= v-k-1+(1+\th)(1+\tau)$.
\item
\label{SRGd}
$x=-1$ and 
\[
y= \frac{-v+2\th^2+2 \pm \sqrt{(v-4)(v-4\th^2)}}{ 2(1+\th)(1+\tau)},
\]
and $A(X)$ is the incidence matrix of a symmetric $(v, k, k+\th\tau)$-design.

\item
\label{SRGe}
$x+x\inv$ is a zero of the quadratic 
\begin{equation}
-\th\tau z^2+\al z+\be+2\th\tau
\label{Eqn_SRGe}
\end{equation}
with 
\begin{align*}
\al &=v(\th +\tau+1)+(\th +\tau)^2,\\
\be &=v+v(1+\th +\tau)^2-2\th ^2-2\th \tau-2\tau^2 \nonumber
\end{align*}
and
\[
y=\frac{\th\tau x^3-[v(\th+\tau+1)-2\th-2\tau-1]x^2-(v+2\th+2\tau+\tau\th)x-1}{(x^2-1)(1+\th)(1+\tau)}.
\]
\end{enumerate}
\qed
\end{theorem}

When $v>4$, $|\frac{1}{2}(2-v\pm\sqrt{v^2-4v})|\neq 1$, and any complex Hadamard matrix in the adjacency algebra
results from Cases~(\ref{SRGc}), (\ref{SRGd}) or (\ref{SRGe}) of Theorem~\ref{Thm_SRG}.

\begin{lemma}
\label{Lem_SRG}
Let $X$ be a primitive strongly regular 
graph with at least five vertices where $v\geq 2k$.
Let $W =I+xA(X)+yA(\overline{X})$
where $(x,y)$ is a solution given by Theorem~\ref{Thm_SRG}.
If $|x|=1$ then 
$\th+\tau \in \{-2, -1, 0\}.$
\end{lemma}
\proof
In Theorem~\ref{Thm_SRG},  $\th+\tau=-2$ in Case~(\ref{SRGc}) and $\th+\tau=0$ in Case~(\ref{SRGd}),
so we need to consider only Case~(\ref{SRGe}).

Let $\Delta = \al^2 +4\th\tau(\be+2\th\tau)$.
Solving the quadratic in (\ref{Eqn_SRGe}) gives
\begin{equation*}
x+x\inv = \frac{-\al \pm \sqrt{\Delta}}{-2\th\tau}.
\end{equation*}
If $|x|=1$ then $x+x\inv\in \reals$ 
and $|x+x\inv |\leq 2$.

We assume $\Delta\geq 0$ for the rest of this proof.

If $\th+\tau \geq 1$, then
\begin{eqnarray*}
\al + 4\th\tau +\sqrt{\Delta}
&=& v(\th+\tau+1)+(\th+\tau)^2  +4\th\tau +\sqrt{\Delta}\\
&\geq& 2v+1+4\th\tau +\sqrt{\Delta}\\
&\geq& 4(k+\th\tau)+1 +\sqrt{\Delta}.
\end{eqnarray*} 
Since $k+\th\tau=c$, we have 
$\al + 4\th\tau +\sqrt{\Delta} >0$.
It follows from $\th\tau <0$ that
\begin{equation*}
\frac{-\al -\sqrt{\Delta}}{-2\th\tau} < -2.
\end{equation*}
From
\begin{equation*}
\left((\al +4\th\tau) - \sqrt{\Delta}\right) 
\left((\al +4\th\tau) + \sqrt{\Delta}\right) 
= -4\th\tau(v-4)(\th+\tau)^2 >0,
\end{equation*}
we have
$\al + 4\th\tau -\sqrt{\Delta} > 0$ and
\begin{equation*}
\frac{-\al +\sqrt{\Delta}}{-2\th\tau} < -2.
\end{equation*}
We conclude that if $\th+\tau\geq 1$
then $x+x\inv <-2$ and $|x|\neq 1$.

Suppose $\th +\tau \leq -3$.
Now
$\al -4\th\tau = (\th +\tau)(v+\th+\tau) +v -4\th\tau$.
Since $v+\th+\tau = v+a-c >0$, we have
\begin{eqnarray*}
\al -4\th\tau -\sqrt{\Delta}
&\leq & -3(v+\th+\tau) +v -4\th\tau.
\end{eqnarray*}
Substituting $\th+\tau=a-c$ and $\th\tau=c-k$ gives
\begin{equation*}
\al -4\th\tau -\sqrt{\Delta} \leq -2(v-2k)-c-3a<0.
\end{equation*}
It follows that
\begin{equation*}
\frac{-\al +\sqrt{\Delta}}{-2\th\tau} > 2.
\end{equation*}
From
\begin{equation*}
\left((\al -4\th\tau) - \sqrt{\Delta}\right) 
\left((\al -4\th\tau) + \sqrt{\Delta}\right) 
= -4\th\tau v(\th+\tau+2)^2 >0,
\end{equation*}
we see that
$ \al -4\th\tau +\sqrt{\Delta} < 0$, so
\begin{equation*}
\frac{-\al -\sqrt{\Delta}}{-2\th\tau} > 2.
\end{equation*}
Therefore $x+x\inv>2$ and $|x|\neq 1$ when $\th+\tau\leq -3$.
\qed

\subsection{Parameters}
In this section, we let $X$ be a primitive strongly regular 
graph with at least five vertices, and
let 
$W =I+xA(X)+yA(\overline{X})$
where $(x,y)$ is a solution given by Theorem~\ref{Thm_SRG}.
\begin{lemma}
\label{Lem_FlatSRG1to3}
Suppose $\th+\tau=0$.  Then $W$ is a complex Hadamard matrix
if and only if $X$ has one of the following parameters:
\begin{enumerate}[i.]
\item
\label{L_SRGFlat1}
$(4\th^2,2\th^2-\th,\th^2-\th,\th^2-\th)$
\item
\label{L_SRGFlat2}
$(4\th^2,2\th^2+\th,\th^2+\th,\th^2+\th)$
\item
\label{L_SRGFlat3}
$(4\th^2-1,2\th^2,\th^2,\th^2)$
\end{enumerate}
\end{lemma}
\proof
Substituting $\tau=-\th$ to
the quadratic in (\ref{Eqn_SRGe}) yields
\begin{equation*}
x+x\inv = -2 \quad \text{or}\quad x+x\inv = \frac{(-v+2\th^2)}{\th^{2}}.
\end{equation*}
The former gives $x=-1$, and from Theorem~\ref{Thm_SRG}~(\ref{SRGd}), we get
\begin{equation*}
y = \frac{-v+2\th^2+2 \pm \sqrt{(v-4)(v-4\th^2)}}{2(1-\th^2)}.
\end{equation*}
Note that $v< 4\th^2$ if and only if
$y \not\in \reals$ and $|y|=1$.

Furthermore, $y=1$ is equivalent to $(\th^2-1)(v-4\th^2)=0$.
If $\th=-\tau=1$ then $v=k+1$ and $X$ is complete.
We conclude that $y=1$ if and only if $v=4\th^2$.
On the other hand, 
$y=-1$ gives $(v-4)(1-\th^2)=0$ which never occurs.
We see that $|y|=1$ if and only if $v\leq 4\th^2$.

When
\begin{equation*}
x+x\inv = \frac{-v+2\th^2}{\th^2},
\end{equation*}
we have
$|x+x\inv| \leq 2$
if and only if $0\leq v\leq 4\th^2$.
Thus $|x|=1$ when $v\leq 4\th^2$.

Simplifying
\begin{equation*}
v=\frac{(k-\th)(k+\th)}{k-\th^2} \leq 4\th^2
\end{equation*}
gives $[k-(2\th^2-\th)][k-(2\th^2+\th)] \leq 0$.
Then we have 
\begin{equation*}
2\th^2-\th \leq k \leq 2\th^2+\th.
\end{equation*}
The integrality of $m_{\th} = \frac{1}{2}(v-1-\frac{k}{\th})$
implies that $\th$ divides  $k$.
Hence $k$ can be $2\th^2-\th$, $2\th^2$ or $2\th^2+\th$,
which give parameter set in (\ref{SRGFlat1}), (\ref{SRGFlat3}) 
and (\ref{SRGFlat2}), respectively.
Now (\ref{SRGFlat1}) and (\ref{SRGFlat2}) are parameter sets for
Latin-square graphs $L_{\th}(2\th)$ and negative Latin-square graphs
$NL_{\th}(2\th)$, respectively.  
Theomem~(\ref{Thm_SRG}) gives only Hadamard matrices,
as expected from \cite{MR0282872}, 
For parameter set (\ref{SRGFlat3}), we get
\begin{equation}
\label{Eqn_SRG1a}
x=-1 \quad \text{and} \quad
y=\frac{2\th^2-3\pm\sqrt{4\th^2-5}i}{2(\th^2-1)},
\end{equation}
or
\begin{equation}
\label{Eqn_SRG1b}
x=\frac{-2\th^2+1\pm\sqrt{4\th^2-1}i}{2\th^2}
\end{equation}
and by Theorem~\ref{Thm_SRG}~(\ref{SRGe})
\begin{equation*}
y=
\frac{1}{x-x\inv}\left(\frac{\th\tau x-1}{(1+\th)(1+\tau)}
(x+x\inv-2+v)-(v-2)x-2\right)=1.
\end{equation*}
In all cases, $|x|=|y|=1$.
\qed

\begin{lemma}
\label{Lem_FlatSRG4to5}
Suppose $\th+\tau=-1$.  Then $W$ is a complex Hadamard matrix
if and only if $X$ or $\overline{X}$ has one of the following parameters:
\begin{enumerate}[i.]
\item
\label{L_SRGFlat4}
$(4\th^2+4\th+1,2\th^2+2\th, \th^2+\th-1,\th^2+\th)$.
\item
\label{L_SRGFlat5}
$(4\th^2+4\th+2,2\th^2+\th,\th^2-1,\th^2)$
\end{enumerate}
\end{lemma}
\proof
If $X$ is a conference graph then we get the parameter set
(\ref{L_SRGFlat4}) and
\begin{equation}
\label{Eqn_SRG2}
x=y\inv=\frac{1\pm\sqrt{(2\th+1)(2\th+3)}i}{2(\th+1)}
\end{equation}
or
\begin{equation}
\label{Eqn_SRG3}
x=y\inv=\frac{-1\pm\sqrt{(2\th+1)(2\th-1)}i}{2\th}.
\end{equation}
Both solutions satisfy $|x|=|y|=1$.

We assume $\th$ is a positive integer.
Theorem~\ref{Thm_SRG}~(\ref{SRGe}) gives
\begin{equation*}
x+x\inv = \frac{-1\pm\sqrt{(2\th+1)^4-4v\th(\th+1)}}{2\th(\th+1)}
\end{equation*}
which is real if and only if
\begin{equation*}
v \leq \frac{(2\th+1)^4}{4\th(\th+1)} =
4\th^2+4\th+2+\frac{1}{4\th(\th+1)}
\end{equation*}
We see that
\begin{equation*}
v=\frac{(k-\th)(k+\th+1)}{k-\th(\th+1)} \leq 4\th^2+4\th+2.
\end{equation*}
After simplification, this inequality becomes 
$(k-2\th^2-\th)(k-2\th^2-3\th-1) \leq 0$.
Hence we have $2\th^2+\th \leq k \leq 2\th^2+3\th+1$.

Now write $k=2\th^2+\th +h$ for some integer $0\leq h \leq 2\th+1$.
Then
\begin{equation*}
v=\frac{(2\th^2+h)(2\th^2+2\th+1+h)}{\th^2+h} = 4\th^2+4\th+2 + \frac{h(h-2\th-1)}{\th^2+h}.
\end{equation*}
For $v \in \ints$, 
either $h=0$, $h=2\th+1$, or 
$(\th^2+h)$ divides $(h(h-2\th-1))$.
If $h(h-2\th-1)\neq 0$ then
the last condition implies $|\th^2+h| \leq |h(h-2\th-1)|$.
Since
\begin{equation*}
|\th^2+h|-|h(h-2\th-1)|=(\th^2+h)+h(h-2\th-1) = (\th-h)^2,
\end{equation*}
the inequality holds if and only if $h=\th$.

When $h=\th$ and $k=2\th^2+2\th$, we get parameter set (\ref{L_SRGFlat4}) and
$X$ is a conference graph.
When $h=2\th+1$, $\overline{X}$ has parameter set (\ref{L_SRGFlat5}).
When $h=0$, $X$ has parameter set (\ref{L_SRGFlat5}).
In the last two cases, Equation~(\ref{Eqn_SRGe}) gives $x+x\inv =0$ or $-\th\inv(1+\th)\inv$.
Hence we have
\begin{equation*}
x=\pm i \quad \text{and} \quad y=x\inv,
\end{equation*}
or
\begin{equation*}
x=\frac{-1\pm\sqrt{4\th^2(\th+1)^2-1}\ i}{2\th(\th+1)}
\quad \text{and} \quad 
y=x\inv=\frac{-1\mp\sqrt{4\th^2(\th+1)^2-1}\ i}{2\th(\th+1)}
\end{equation*}
In all cases, $|x|=|y|=1$.
\qed

\begin{theorem}
\label{Thm_FlatSRG}
Suppose $X$ is a primitive strongly regular graph with at least five vertices.
Let $k$, $\th$ and $\tau$ be the eigenvalues of $A(X)$ satisfying $k>\th>0>\tau$.
Let 
\begin{equation*}
W =I+xA(X)+yA(\overline{X})
\end{equation*} 
where $(x,y)$ is a solution from
Theorem~\ref{Thm_SRG}. 

Then $W$ is a complex Hadamard matrix if and only if $X$ or $\overline{X}$
has one of the following parameters:
\begin{enumerate}[i.]
\item
\label{SRGFlat1}
$(4\th^2,2\th^2-\th,\th^2-\th,\th^2-\th)$
\item
\label{SRGFlat2}
$(4\th^2,2\th^2+\th,\th^2+\th,\th^2+\th)$
\item
\label{SRGFlat3}
$(4\th^2-1,2\th^2,\th^2,\th^2)$
\item
\label{SRGFlat4}
$(4\th^2+4\th+1,2\th^2+2\th, \th^2+\th-1,\th^2+\th)$.
\item
\label{SRGFlat5}
$(4\th^2+4\th+2,2\th^2+\th,\th^2-1,\th^2)$
\end{enumerate}
\end{theorem}
\proof
Observe that $\th+\tau = 0$ if and only if 
$(-1-\tau)+(-1-\th)=-2$.  
By Lemma~\ref{Lem_SRG}, it is sufficient to 
consider only the strongly regular graphs with $\th+\tau=0$ or $\th+\tau=-1$.
The theorem follows from Lemma~\ref{Lem_FlatSRG1to3} and Lemma~\ref{Lem_FlatSRG4to5}.
\qed

\subsection{Remarks}
We remark on the parameter sets (\ref{SRGFlat3}) to (\ref{SRGFlat5})
in Theorem~\ref{Thm_FlatSRG},
which give complex Hadamard matrices with non-real entries.
\begin{enumerate}
\item
Let $X$ be a strongly regular graph with
parameters $(4\th^2-1,2\th^2,\th^2,\th^2)$.

Let $S(X)=J-I-2A(X)$ be the Seidel matrix of $X$.
Then
\begin{equation*}
T = 
\begin{pmatrix}
0 & \bf{1}^T\\
\bf{1}&S(X)
\end{pmatrix}
\end{equation*}
is a regular two-graph with minimal polynomial $z^2+2z-(4\th^2-1)$.
It is a generalized conference matrix and $T+I$ is
a Hadamard matrix.
Solution~(\ref{Eqn_SRG1a}) gives a complex Hadamard matrix that
is equivalent to the construction from $T+I$
described in Theorem~4.1 of \cite{Exotic}.

Note that $A(X)$ is the incidence matrix of the complement 
of a Hadamard design.
Solution~(\ref{Eqn_SRG1b}) gives the same construction as in Theorem~2.4
of \cite{Exotic} and Theorem~6.1 of \cite{TypeIIComb}.

The complex Hadamard matrices obtained from these graphs, for $\th \geq 2$,
do not result from \Dita's construction, see Proposition~\ref{Prop_SRGiii} of Section~\ref{Nomura}.

\item
Let $X$ be a conference graph with
parameters 
\begin{equation*}
(4\th^2+4\th+1,2\th^2+2\th, \th^2+\th-1,\th^2+\th).
\end{equation*}

Let $S(X)$ be the Seidel matrix of $X$.
Then
\begin{equation*}
T = 
\begin{pmatrix}
0 & \bf{1}^T\\
\bf{1}&S(X)
\end{pmatrix}
\end{equation*}
is a regular two-graph with minimal polynomial $z^2-(2\th+1)^2$.
Thus $T$ is a symmetric conference matrix and the solutions in (\ref{Eqn_SRG2})
and (\ref{Eqn_SRG3})
are the same as the one given by Theorem~3.1 of \cite{Exotic}.

When $\th^2+\th > 2$,
the complex Hadamard matrices obtained from these graphs 
do not come from \Dita's construction, see Theorem~3.2 of \cite{MR2106965}.

When $\th^2+\th=2$, Solution~(\ref{Eqn_SRG3}) gives the Kronecker product
of two Potts model of size three.  This is the only complex Hadamard matrix
obtained from these graphs that results from \Dita's construction.

\item
Let $X$ be a strongly regular graph with parameter set
\begin{equation*}
(4\th^2+4\th+2,2\th^2+\th,\th^2-1,\th^2).
\end{equation*}

The Seidel matrix $S(X)$ of $X$ is a regular two-graph with minimal polynomial
$z^2-(2\th+1)^2=0$.
Hence $S(X)$ is a regular symmetric conference matrix.
The complex Hadamard matrices with $x=y\inv=\pm i$
were constructed by Turyn in \cite{MR0270938}.

Conversely, every regular conference matrix is symmetric \cite{MR1284019}
and it is the Seidel matrix
of a strongly regular graph $X$ where $X$ or $\overline{X}$ has the above
parameters.
\begin{corollary}
Given a regular conference matrix $C$ of size $4\th^2+4\th+2$.
Let $W$ be the matrix obtained from $C$ by replacing $0$, $1$ and $-1$
by 
\begin{equation*}
1, \quad
x=\frac{-1\pm\sqrt{4\th^2(\th+1)^2-1}\ i}{2\th(\th+1)},
\quad \text{and} \quad 
\overline{x}
\end{equation*}
respectively.  Then $W$ is a complex Hadamard matrix.
\end{corollary}
None of these complex Hadamard matrices results from \Dita's construction,
see Proposition~\ref{Prop_SRGv} of Section~\ref{Nomura}.
\end{enumerate}

\section{Covers of Complete Graphs}
\label{cover}
A graph $X$ of diameter $d$ is {\sl antipodal} if the relation 
"at distance $0$ or $d$" is an equivalence relation.
We are interested in the distance regular graphs with diameter three
that are antipodal.  These are the distance regular covers of the complete
graph.
There are three parameters $(n,r,c_2)$ associated with each of these graphs
where $n$ is the number of equivalence classes;
$r$ is the size of each equivalence class;
and every pair of vertices at distance two have exactly $c_2$ common neighbours.
The intersection numbers of $X$ depend on $n$, $r$ and $c_2$,
in particular, adjacent vertices have exactly $a_1=n-2-(r-1)c_2$
common neighbours.

Let $X$ be an antipodal distance-regular graph of diameter three 
of parameters $(n,r,c_2)$.
The adjacency matrix of $X$ has four eigenvalues:
\begin{equation*}
n-1,\quad -1,\quad \th = \frac{\delta+\sqrt{\delta^2+4(n-1)}}{2},
\quad\text{and}\quad
\tau = \frac{\delta-\sqrt{\delta^2+4(n-1)}}{2},
\end{equation*}
where $\delta = n-2-rc_2$.

See \cite{MR1186756} for details on distance regular covers of
the complete graph including the next result.
\begin{lemma} 
If $\delta=0$ then $\th=-\tau = \sqrt{n-1}$.
Otherwise, $\th$ and $\tau$ are integers.
\qed
\end{lemma}

From Theorem~9.1 and its proof in \cite{TypeIIComb}, we see that
\begin{lemma}
\label{Lem_CoverSys}
Suppose $X$ is  
an antipodal distance regular graph of diameter three with parameters $(n,r,c_2)$.
For $i=1,2,3$, let $A_i$ be the $i$-th distance matrix of $X$.
Then
\begin{equation*}
W=I+xA_1+yA_2+zA_3
\end{equation*}
is type~II if and only if the following system of equations have a solution:
\begin{align}
(1-x-(r-1)y+(r-1)z) 
\left(1-\frac{1}{x}-\frac{(r-1)}{y}+\frac{(r-1)}{z}\right) &= nr, \label{Eqn_Cover1}\\
(1+\th x-\th y-z)(1+\frac{\th}{x}-\frac{\th}{y}-\frac{1}{z}) &=nr, \label{Eqn_Cover2}\\
(1+\tau x-\tau y-z) (1+\frac{\tau}{x}-\frac{\tau}{ y}-\frac{1}{z}) &=nr. \label{Eqn_Cover3}
\end{align}
\qed
\end{lemma}
For complex Hadamard matrices, we add the condition $|x|=|y|=|z|=1$.
\begin{lemma}
\label{Lem_Cover}
If the system of equations in Lemma~\ref{Lem_CoverSys} 
has a solution satisfying
$|x|=|y|=|z|=1$ then $n\leq 16$.
\end{lemma}
\proof
Assume $|x|=|y|=|z|=1$. 
Then 
for $\alpha \in \{x,y,z,\frac{x}{y},\frac{x}{z},\frac{y}{z}\}$,
we have $\al+\al\inv \in \reals$ and $|\al +\al\inv|\leq 2$.

Expanding the right-hand side of Equation~(\ref{Eqn_Cover1}) gives
\begin{eqnarray}
\nonumber
nr&=&2+2(r-1)^2-(x+\frac{1}{x})-(r-1)(y+\frac{1}{y})+(r-1)(z+\frac{1}{z})\\
&&+(r-1)(\frac{y}{x}+\frac{x}{y})
\nonumber
- (r-1)(\frac{x}{z}+\frac{z}{x})- (r-1)^2(\frac{z}{y}+\frac{y}{z})\\
&\leq & 4r^2.
\label{Eqn_FlatCover1}
\end{eqnarray}
Similarly, Equations~(\ref{Eqn_Cover2}) and (\ref{Eqn_Cover3}) yield
\begin{eqnarray}
\nonumber
nr &=& 2+2\th^2+\th(x+\frac{1}{x})-\th(y+\frac{1}{y})-(z+\frac{1}{z})
-\th^2(\frac{y}{x}+\frac{x}{y}) -\th(\frac{z}{x}+\frac{x}{z})+
\th(\frac{z}{y}+\frac{y}{z})\\
&\leq& 4(1+\th)^2
\label{Eqn_FlatCover2}
\end{eqnarray}
and
\begin{eqnarray}
\nonumber
nr &=& 2+2\tau^2+\tau(x+\frac{1}{x})-\tau(y+\frac{1}{y})-(z+\frac{1}{z})
 -\tau^2(\frac{y}{x}+\frac{x}{y}) -\tau(\frac{z}{x}+\frac{x}{z})+
\tau(\frac{z}{y}+\frac{y}{z})\\
\nonumber
&\leq& 2+2\tau^2-2\tau-2\tau +2+2\tau^2-2\tau-2\tau\\
&=& 4(1-\tau)^2,
\label{Eqn_FlatCover3}
\end{eqnarray}
respectively.

When $\delta \geq 2$, the expression
\begin{eqnarray*}
n^2-16(1-\tau)^2 
&=&n^2-4(\delta-2)^2-4(\delta^2+4(n-1))+
8(\delta-2)\sqrt{\delta^2+4(n-1)}\\
&\geq &n^2-4(\delta-2)^2-4(\delta^2+4(n-1))+ 8(\delta-2)\delta\\
&=& n(n-16).
\end{eqnarray*}
If $n\geq 17$ then $n^2 > 16(1-\tau)^2$,
it follows from (\ref{Eqn_FlatCover3}) that $n > 4r$
which contradicts (\ref{Eqn_FlatCover1}).

When $\delta=1$, we have $\tau=\frac{1-\sqrt{4n-3}}{2}$ and
\begin{equation*}
nr \leq 4(1-\tau)^2 =4n-2+2\sqrt{4n-3} < 4n+4\sqrt{n}.
\end{equation*}
When $n\geq 17$, $r \leq \lfloor 4+\frac{4}{\sqrt{n}}\rfloor=4$
and (\ref{Eqn_FlatCover1}) does not hold.

When $\delta=0$, $\th=-\tau=\sqrt{n-1}$.
Adding both sides of Equations~(\ref{Eqn_Cover2}) and (\ref{Eqn_Cover3}) gives
\begin{equation*}
2nr = 4n-2(z+\frac{1}{z})-2(n-1)(\frac{y}{x}+\frac{x}{y})\\
\leq 8n.
\end{equation*}
Hence $r \leq 4$ and (\ref{Eqn_FlatCover1}) does not hold for $n\geq 17$.

When $\delta=-1$, we have $\th=\frac{-1+\sqrt{4n-3}}{2}$ and
\begin{equation*}
nr \leq 4(1+\th)^2 =4n-2+2\sqrt{4n-3}.
\end{equation*}
Using the same argument for $\delta=1$, we see that if $n\geq 17$ then 
$nr \leq 4(1+\th)^2$ implies $n>4r$.

When $\delta \leq -2$, the expression
\begin{eqnarray*}
n^2-16(1+\th)^2
&=&
n^2-8\delta^2-16\delta-16n-8(\delta+2)\sqrt{\delta^2+4(n-1)}\\
&\geq & n^2-16n-8\delta^2-16\delta-8(\delta+2)(-\delta)\\
&=&n(n-16).
\end{eqnarray*}
If $n\geq 17$ then $n^2 > 16(1+\th)^2$,
but (\ref{Eqn_FlatCover2}) gives $n > 4r$.

We conclude that when $n\geq 17$, 
(\ref{Eqn_FlatCover1}),
(\ref{Eqn_FlatCover2}) and (\ref{Eqn_FlatCover3}) do not hold simultaneously with $|x|=|y|=|z|=1$.
\qed

\begin{theorem}
\label{Thm_Cover}
There are only finitely many antipodal distance-regular graphs of
diameter three where the span of $\{I, A_1, A_2, A_3\}$ contains a complex
Hadamard matrix.
\end{theorem}
\proof
Since $a_1 \geq 0$, we have $(r-1)c_2 \leq n-2$.  
Therefore, there are only finitely many parameter sets $(n, r, c_2)$  where 
$n\leq 16$.
\qed

Computations in Maple revealed only three graphs that give complex
Hadamard matrices:
\begin{enumerate}
\item
The cycle on six vertices: its parameters are $(3,2,1)$ and
\begin{equation*}
y = \frac{-1\pm \sqrt{3}\ i}{2}, \quad z=\pm i, \quad x=yz.
\end{equation*}
We get the matrix $F_6$ in \cite{MR2244963}.
\item
The cube: its parameters are $(4,2,2)$ and
$x=\pm i$, $y=-1$ and $z=-x$.
We get the complex Hadamard matrix
\begin{equation*}
\begin{pmatrix}
1 & i \\i & 1
\end{pmatrix}^{\otimes 3}
\end{equation*}
\item
The line graph of the Petersen graph:
its parameters are $(5,3,1)$ and
\begin{itemize}
\item
$x=1$, $y = \frac{-7\pm \sqrt{15}\ i}{8} $, $z=1$;
\item
$x=\frac{5\pm\sqrt{11}\ i}{6}$, $y = -1$, $z=x$;
\item
$x=\frac{-1\pm\sqrt{15}\ i}{4}$, $y = x\inv$, $z=1$; or
\item
$x=\frac{(3+\sqrt{201})z+15+5\sqrt{201}}{96}$, $y=\frac{-(15+5\sqrt{201})z-3-\sqrt{201}}{96}$, 
$z=\frac{53-3\sqrt{201}\pm \sqrt{-4218+318\sqrt{201}} \ i}{20}$.
\end{itemize}
The first two solutions give complex Hadamard matrices
that are equivalent to $U_{15}$ and $V_{15}$ in \cite{Exotic}.
The last solution is a special case  of Theorem 2 (vi)  in \cite{MR3349366} with $q=4$.
Computation in Maple showed that 
none of these matrices results from \Dita's construction.
\end{enumerate}

\section{Nomura Algebra}
\label{Nomura}

In this section, we consider a complex Hadamard matrix $W$ in the adjacency algebra of a strongly regular graph with parameters
given in Theorem~\ref{Thm_FlatSRG}~(\ref{SRGc}) to (\ref{SRGe}).   We show that the Nomura algebra of $W$ has dimension $2$,
which implies that $W$ cannot be a result of \Dita's construction.

Given an $v\times v$ type II matrix$W$,  define a column vector $\y{\alpha}{\beta}$ as
\begin{equation}
\label{Eqn_Yab}
\y{\alpha}{\beta}(x) = \frac{W(x,\alpha)}{W(x,\beta)}, \quad \text{for $x=1,\ldots,v$,}
\end{equation}
for $\alpha, \beta = 1, \ldots, v$.

The {\sl Nomura algebra of $W$} is the set
\begin{equation*}
\nomw = \{M \ |\ \text{$\y{\alpha}{\beta}$ is an eigenvector of $M$, for $\alpha, \beta =1,\ldots,v$}\}.
\end{equation*}
We see immediately that $I \in \nomw$.  It follows from the type II condition of $W$ that $\nomw$ also contains $J$,
hence $\dim \nomw \geq 2$.   Jaeger, Matsumoto and Nomura proved the $\nomw$ is the Bose-Mesner algebra of
an association scheme \cite{MR1635553}.

Given $n\times n$ matrices $U_1,\ldots, U_m$ and $m\times m$ matrices $V_1,\ldots, V_n$.
The {\sl generalized tensor product of $U_1, \ldots, U_m$ and $V_1, \ldots, V_n$},
denoted by $(U_1,\ldots,U_m) \otimes (V_1,\ldots,V_n)$,
is the $mn\times mn$ matrix whose $(i,j)$-block is
\begin{equation*}
\Delta_{i,j}V_j
\end{equation*}
where
\begin{equation*}
\Delta_{i,j}(h,k)=
\begin{cases}
U_h(i,j) & \text{if $h=k$,}\\
0 & \text{otherwise.}
\end{cases}
\end{equation*}

\begin{lemma} \cite{MR1958007}
\label{Lem_HS}
Let $U_1, \ldots, U_m$ be $n\times n$ matrices and
$V_1,\ldots,V_n$ be $m\times m$ matrices.
Then $(U_1,\ldots,U_m) \otimes (V_1,\ldots,V_n)$ is a type II matrix if
and only if
$U_1, \ldots, U_m, V_1,\ldots, V_n$ are type II matrices.
\end{lemma}
Here we use $I_m$ to denote the $m\times m$ identity matrix and $J_n$ to denote
the $n\times n$ matrix of all ones.
\begin{theorem} \cite{MR1958007}
\label{Thm_HS}
Let $W$ be a $mn \times mn$ type II matrix. Then the following are equivalent.
\begin{enumerate}
\item
There exists a permutation matrix $P$ such that $P(J_n\otimes I_m)P^T \in \nomw$.
\item
$W$ is type II equivalent to a generalized tensor product 
$(U_1,\ldots,U_m) \otimes (V_1,\ldots,V_n)$.
\end{enumerate}
\end{theorem}
If $U_1,\ldots, U_m$ 
and $V_1,\ldots, V_n$ are complex Hadamard matrices of order $n$ and $m$, respectively,
then $(U_1,\ldots,U_m) \otimes (V_1,\ldots,V_n)$ a complex Hadamard matrix
of order $mn$.

Observe that when $U_1=\ldots=U_m$, we have
\begin{equation*}
\Delta_{i,j} = U_1(i,j) I_m,
\end{equation*}
and the generalized tensor product $(U_1,\ldots,U_m) \otimes (V_1,\ldots,V_n)$
is a complex Hadamard matrix of Di{\c{t}}{\u{a}}-type \cite{MR2065675}.
\begin{corollary}
Let $W$ be a complex Hadamard matrix.  If $W$ is of Di{\c{t}}{\u{a}}-type then
$\nomw$ has dimension at least three. 
\end{corollary}

We adapt the computation in \cite{MR2106965} to the strongly regular graphs in Theorem~\ref{Thm_FlatSRG}~(\ref{SRGc}) to (\ref{SRGe}).
The following is Lemma~3.1 of \cite{MR2106965} which was phrased for conference
graphs but the statement applies to all strongly regular graphs.
\begin{lemma}
\label{Lem_Trivial}
Let $X$ be a strongly regular graph on $v\geq 5$ vertices with degree $k\geq 2$.
Let $W=I+xA(X)+yA(\overline{X})$ be a complex Hadamard matrix.
If the Hermitian product $<\y{\a}{\b},\y{\g}{\a}>$ is non-zero for all
adjacent vertices $\a$ and $\b$ and for all $\g\neq \a, \b$, then
\begin{equation*}
\nomw = \spn \{I, J\}.
\end{equation*}
\qed
\end{lemma}

Let $\a,\b$ and $\g$ be vertices of a strongly regular graph $X$ with parameters 
$(v,k,a,c)$. 
For $x_{\a} \in \{\a,\ba\}$, $x_{\b} \in \{\b,\bb\}$ and $x_{\g} \in \{\g,\bg\}$, 
we define $\C{x_{\a}}{x_{\b}}{x_{\g}}$ to be the set of vertices  
that are adjacent (not adjacent) to $x_v$  
if $x_v=v$ ($x_v=\bar{v}$, respectively), for $v=\a,\b,\g$. 
For instance $\C{\a}{\b}{\g}$ is the set of common neighbours of $\a$, $\b$ and  $\g$ 
while $\C{\a}{\b}{\bg}$ is the set of common neighbours of $\a$ and $\b$ that  are not adjacent to $\g$. 
Now the vertex set of $X$ is partitioned into 
\begin{equation*} 
\{\a,\b,\g\} \cup \C{\a}{\b}{\g} \cup \C{\ba}{\b}{\g} \cup 
\C{\a}{\bb}{\g} \cup \C{\a}{\b}{\bg} \cup \C{\a}{\bb}{\bg} 
\cup \C{\ba}{\b}{\bg} \cup \C{\ba}{\bb}{\g} \cup \C{\ba}{\bb}{\bg}. 
\end{equation*} 
Let $W=I+xA(X)+yA(\overline{X})$ be a complex Hadamard matrix.
From (\ref{Eqn_Yab}), we have 
\begin{equation*} 
\y{\a}{\b}(u)\overline{\y{\g}{\a}(u)} = \y{\a}{\b}(u) \y{\a}{\g}(u)
= \frac{W(u,\a)^2}{W(u,\b)W(u,\g)}. 
\end{equation*} 
It is easy to verify that  
\begin{equation*} 
\y{\a}{\b}(u)\overline{\y{\g}{\a}(u)} = 
\begin{cases} 
1      & \text{if } u \in \C{\a}{\b}{\g} \cup \C{\ba}{\bb}{\bg}\\ 
xy\inv    & \text{if } u \in \C{\a}{\b}{\bg} \cup \C{\a}{\bb}{\g}\\ 
x\inv y& \text{if } u \in \C{\ba}{\b}{\bg} \cup \C{\ba}{\bb}{\g}\\ 
x^{-2}y^2 & \text{if } u \in \C{\ba}{\b}{\g}\\ 
x^2y^{-2}    & \text{if } u \in \C{\a}{\bb}{\bg}.\\ 
\end{cases} 
\end{equation*} 
Hence the Hermitian product 
\begin{eqnarray} 
\label{Eqn_Herm} 
<\y{\a}{\b},\y{\g}{\a}>
&=& \frac{1}{W(\a,\b)W(\a,\g)} + 
\frac{W(\b,\a)^2}{W(\b,\g)} + \\ 
\nonumber 
&&\frac{W(\g,\a)^2}{W(\g,\b)} +  
|\C{\a}{\b}{\g} \cup \C{\ba}{\bb}{\bg}|+  
|\C{\a}{\bb}{\g} \cup \C{\a}{\b}{\bg}| xy\inv +\\ 
\nonumber 
&&|\C{\ba}{\b}{\bg} \cup \C{\ba}{\bb}{\g}| x\inv y + 
|\C{\ba}{\b}{\g}| x^{-2}y^2 + |\C{\a}{\bb}{\bg}| x^2y^{-2}.  
\end{eqnarray}  
In the following computation, we let $\a$ and $\b$ be adjacent vertices in 
$X$.  
\begin{enumerate}
\item
We first consider the case where $\g$ is adjacent to both $\a$ and $\b$.  
Then $W(\a,\b)=W(\a,\g)=W(\b,\g)=x$.  
We use $\N{v}$ to denote the set of neighbours of $v$ in $X$. 
Then we get 
\begin{footnotesize}
\begin{eqnarray*} 
\N{\a} &=& 
\C{\a}{\b}{\g} \cup \C{\a}{\bb}{\g} \cup \C{\a}{\b}{\bg} \cup 
\C{\a}{\bb}{\bg} \cup \{ \b,\g\} \\ 
\N{\b} &=& 
\C{\a}{\b}{\g} \cup \C{\ba}{\b}{\g} \cup \C{\a}{\b}{\bg} \cup 
\C{\ba}{\b}{\bg} \cup \{ \a,\g\}\\ 
\N{\g} &=& 
\C{\a}{\b}{\g} \cup \C{\ba}{\b}{\g} \cup \C{\a}{\bb}{\g} \cup 
\C{\ba}{\bb}{\g} \cup \{ \a,\b\}\\ 
\N{\a} \cap \N{\b} &=& 
\C{\a}{\b}{\g} \cup \C{\a}{\b}{\bg} \cup \{\g\}\\ 
\N{\a} \cap \N{\g} &=& 
\C{\a}{\b}{\g} \cup \C{\a}{\bb}{\g} \cup \{\b\}\\ 
\N{\b} \cap \N{\g} &=& 
\C{\a}{\b}{\g} \cup \C{\ba}{\b}{\g} \cup \{\a\}\\ 
V(X) &=& \C{\a}{\b}{\g} \cup \C{\ba}{\b}{\g} \cup 
\C{\a}{\bb}{\g} \cup \C{\a}{\b}{\bg} \cup \C{\a}{\bb}{\bg} 
\cup \C{\ba}{\b}{\bg} \cup \C{\ba}{\bb}{\g} \cup \C{\ba}{\bb}{\bg} 
\cup \{\a,\b,\g\}. 
\end{eqnarray*} 
\end{footnotesize}
Now we translate the above to the following system of 
equations. 
\begin{footnotesize}
\begin{alignat*}{6}
k &=& 
|\C{\a}{\b}{\g}|& &+ |\C{\a}{\bb}{\g}|& + |\C{\a}{\b}{\bg}|& + 
|\C{\a}{\bb}{\bg}| &&&&+ 2 \\ 
k &=& 
|\C{\a}{\b}{\g}| &+ |\C{\ba}{\b}{\g}|&& + |\C{\a}{\b}{\bg}| &&+ 
|\C{\ba}{\b}{\bg}|&&& + 2 \\ 
k &=& 
|\C{\a}{\b}{\g}| &+ |\C{\ba}{\b}{\g}|& + |\C{\a}{\bb}{\g}| &&&&+ 
|\C{\ba}{\bb}{\g}|&& + 2  \\ 
a &=& 
|\C{\a}{\b}{\g}| &&&+ |\C{\a}{\b}{\bg}| &&&&&+ 1  \\ 
a &=& 
|\C{\a}{\b}{\g}| &&+ |\C{\a}{\bb}{\g}|&&&&&& + 1 \\ 
a &=& 
|\C{\a}{\b}{\g}| &+ |\C{\ba}{\b}{\g}|&&&&&&& + 1 \\ 
v &=& |\C{\a}{\b}{\g}| &+ |\C{\ba}{\b}{\g}|& + 
|\C{\a}{\bb}{\g}|& + |\C{\a}{\b}{\bg}|& + |\C{\a}{\bb}{\bg}|& + 
|\C{\ba}{\b}{\bg}|& + |\C{\ba}{\bb}{\g}|& + |\C{\ba}{\bb}{\bg}|& 
+ 3. 
\end{alignat*} 
\end{footnotesize}
Solving this system of equations, we get 
\begin{eqnarray*} 
|\C{\a}{\b}{\g}|&=&m, \\
|\C{\ba}{\b}{\g}| &=& |\C{\a}{\bb}{\g}|=|\C{\a}{\b}{\bg}|= a-1-m,\\
|\C{\a}{\bb}{\bg}|&=& |\C{\ba}{\b}{\bg}|=|\C{\ba}{\bb}{\g}|=k-2a+m,
\quad \text{and} \quad \\ 
|\C{\ba}{\bb}{\bg}| &=& v-3k+3a-m,
\end{eqnarray*} 
for some integer $m$. 
Using Equation~(\ref{Eqn_Herm}), we have 
\begin{eqnarray} 
\label{Eqn_a_b}
&&<\y{\a}{\b},\y{\g}{\a}> \\
&=&x^{-2}+2x+ (v-3k+3a)+2(a-1-m)xy\inv
 +2(k-2a+m)x\inv y \nonumber\\
&&+ (a-1-m)x^{-2}y^2+(k-2a+m)x^2y^{-2}.
\nonumber
\end{eqnarray} 
\item
Suppose $\g$ is adjacent to $\a$ but not to $\b$
which gives 
\begin{equation*}
W(\a,\b)=W(\a,\g)=x \quad \text{and}\quad W(\b,\g)=y.  
\end{equation*}
We get
\begin{footnotesize}
\begin{alignat*}{6}
k &=& 
|\C{\a}{\b}{\g}|& &+ |\C{\a}{\bb}{\g}|& + |\C{\a}{\b}{\bg}|& + 
|\C{\a}{\bb}{\bg}| &&&&+ 2 \\ 
k &=& 
|\C{\a}{\b}{\g}| &+ |\C{\ba}{\b}{\g}|&& + |\C{\a}{\b}{\bg}| &&+ 
|\C{\ba}{\b}{\bg}|&&& + 1 \\ 
k &=& 
|\C{\a}{\b}{\g}| &+ |\C{\ba}{\b}{\g}|& + |\C{\a}{\bb}{\g}| &&&&+ 
|\C{\ba}{\bb}{\g}|&& + 1  \\ 
a &=& 
|\C{\a}{\b}{\g}| &&&+ |\C{\a}{\b}{\bg}| &&&&&  \\ 
a &=& 
|\C{\a}{\b}{\g}| &&+ |\C{\a}{\bb}{\g}|&&&&&&  \\ 
c &=& 
|\C{\a}{\b}{\g}| &+ |\C{\ba}{\b}{\g}|&&&&&&& + 1 \\ 
v &=& |\C{\a}{\b}{\g}| &+ |\C{\ba}{\b}{\g}|& + 
|\C{\a}{\bb}{\g}|& + |\C{\a}{\b}{\bg}|& + |\C{\a}{\bb}{\bg}|& + 
|\C{\ba}{\b}{\bg}|& + |\C{\ba}{\bb}{\g}|& + |\C{\ba}{\bb}{\bg}|& 
+ 3. 
\end{alignat*} 
\end{footnotesize}
Solving this system of equations yields
\begin{eqnarray*} 
|\C{\a}{\b}{\g}|&=&m, \\
|\C{\ba}{\b}{\g}| &=& c-1-m\\
|\C{\a}{\bb}{\g}|&=&|\C{\a}{\b}{\bg}|= a-m,\\
|\C{\a}{\bb}{\bg}|&=& k-2a-2+m\\
|\C{\ba}{\b}{\bg}|&=&|\C{\ba}{\bb}{\g}|=k-a-c+m,
\quad \text{and} \quad \\ 
|\C{\ba}{\bb}{\bg}| &=& v-3k+2a+c-m,
\end{eqnarray*} 
for some integer $m$. 
Using Equation~(\ref{Eqn_Herm}), we have 
\begin{eqnarray} 
\label{Eqn_a_bb}
&&<\y{\a}{\b},\y{\g}{\a}> \\
&=&
x^{-2}+2x^2y\inv+
 (v-3k+2a+c)+2(a-m)xy\inv\nonumber\\
&& +2(k-a-c+m)x\inv y 
+ (c-1-m)x^{-2}y^2+(k-2a-2+m)x^2y^{-2}.
\nonumber
\end{eqnarray} 
\item
Suppose $\g$ is adjacent to $\b$ but not to $\a$
which gives 
\begin{equation*}
W(\a,\b)=W(\b,\g)=x \quad \text{and}\quad W(\a,\g)=y.  
\end{equation*}
We get
\begin{footnotesize}
\begin{alignat*}{6}
k &=& 
|\C{\a}{\b}{\g}|& &+ |\C{\a}{\bb}{\g}|& + |\C{\a}{\b}{\bg}|& + 
|\C{\a}{\bb}{\bg}| &&&&+ 1 \\ 
k &=& 
|\C{\a}{\b}{\g}| &+ |\C{\ba}{\b}{\g}|&& + |\C{\a}{\b}{\bg}| &&+ 
|\C{\ba}{\b}{\bg}|&&& + 2 \\ 
k &=& 
|\C{\a}{\b}{\g}| &+ |\C{\ba}{\b}{\g}|& + |\C{\a}{\bb}{\g}| &&&&+ 
|\C{\ba}{\bb}{\g}|&& + 1  \\ 
a &=& 
|\C{\a}{\b}{\g}| &&&+ |\C{\a}{\b}{\bg}| &&&&&  \\ 
c &=& 
|\C{\a}{\b}{\g}| &&+ |\C{\a}{\bb}{\g}|&&&&&& +1 \\ 
a &=& 
|\C{\a}{\b}{\g}| &+ |\C{\ba}{\b}{\g}|&&&&&&&  \\ 
v &=& |\C{\a}{\b}{\g}| &+ |\C{\ba}{\b}{\g}|& + 
|\C{\a}{\bb}{\g}|& + |\C{\a}{\b}{\bg}|& + |\C{\a}{\bb}{\bg}|& + 
|\C{\ba}{\b}{\bg}|& + |\C{\ba}{\bb}{\g}|& + |\C{\ba}{\bb}{\bg}|& 
+ 3. 
\end{alignat*} 
\end{footnotesize}
Solving this system of equations gives
\begin{eqnarray*} 
|\C{\a}{\b}{\g}|&=&m, \\
|\C{\ba}{\b}{\g}|&=&|\C{\a}{\b}{\bg}|= a-m,\\
|\C{\a}{\bb}{\g}| &=& c-1-m,\\
|\C{\a}{\bb}{\bg}|&=&|\C{\ba}{\bb}{\g}|=k-a-c+m,\\
|\C{\ba}{\b}{\bg}|&=& k-2a-2+m,
\quad \text{and} \quad \\ 
|\C{\ba}{\bb}{\bg}| &=& v-3k+2a+c-m,
\end{eqnarray*} 
for some integer $m$. 
Using Equation~(\ref{Eqn_Herm}), we have 
\begin{eqnarray} 
\label{Eqn_ba_b}
&& <\y{\a}{\b},\y{\g}{\a}> \\
&=& 
x\inv y\inv+x+x\inv y^2 +
(v-3k+2a+c)+(a+c-1-2m)xy\inv\nonumber\\
&& +(2k-3a-c-2+2m)x\inv y 
+ (a-m)x^{-2}y^2+(k-a-c+m)x^2y^{-2} .
\nonumber
\end{eqnarray} 
\item
Suppose $\g$ is not adjacent to $\a$ nor to $\b$
which gives $W(\a,\b)=x$ and $W(\a,\g)=W(\b,\g)=y$.  
We get
\begin{footnotesize}
\begin{alignat*}{6}
k &=& 
|\C{\a}{\b}{\g}|& &+ |\C{\a}{\bb}{\g}|& + |\C{\a}{\b}{\bg}|& + 
|\C{\a}{\bb}{\bg}| &&&&+ 1 \\ 
k &=& 
|\C{\a}{\b}{\g}| &+ |\C{\ba}{\b}{\g}|&& + |\C{\a}{\b}{\bg}| &&+ 
|\C{\ba}{\b}{\bg}|&&& + 1 \\ 
k &=& 
|\C{\a}{\b}{\g}| &+ |\C{\ba}{\b}{\g}|& + |\C{\a}{\bb}{\g}| &&&&+ 
|\C{\ba}{\bb}{\g}|&&   \\ 
a &=& 
|\C{\a}{\b}{\g}| &&&+ |\C{\a}{\b}{\bg}| &&&&&  \\ 
c &=& 
|\C{\a}{\b}{\g}| &&+ |\C{\a}{\bb}{\g}|&&&&&&  \\ 
c &=& 
|\C{\a}{\b}{\g}| &+ |\C{\ba}{\b}{\g}|&&&&&&&  \\ 
v &=& |\C{\a}{\b}{\g}| &+ |\C{\ba}{\b}{\g}|& + 
|\C{\a}{\bb}{\g}|& + |\C{\a}{\b}{\bg}|& + |\C{\a}{\bb}{\bg}|& + 
|\C{\ba}{\b}{\bg}|& + |\C{\ba}{\bb}{\g}|& + |\C{\ba}{\bb}{\bg}|& 
+ 3. 
\end{alignat*} 
\end{footnotesize}
Solving this system of equations, we get 
\begin{eqnarray*} 
|\C{\a}{\b}{\g}|&=&m, \\
|\C{\ba}{\b}{\g}|&=&|\C{\a}{\bb}{\g}|= c-m,\\
|\C{\a}{\b}{\bg}| &=& a-m,\\
|\C{\a}{\bb}{\bg}|&=&|\C{\ba}{\b}{\bg}|=k-a-c-1+m,\\
|\C{\ba}{\bb}{\g}|&=& k-2c+m,
\quad \text{and} \quad \\ 
|\C{\ba}{\bb}{\bg}| &=& v-3k+a+2c-1-m,
\end{eqnarray*} 
for some integer $m$. 
Using Equation~(\ref{Eqn_Herm}), we have 
\begin{eqnarray} 
\label{Eqn_ba_bb}
&&<\y{\a}{\b},\y{\g}{\a}>\\ 
&=& 
x\inv y\inv+y+x^2 y\inv + 
(v-3k+2a+c-1)+(a+c-2m)xy\inv  \nonumber \\
&&+(2k-a-3c-1+2m)x\inv y
+ (c-m)x^{-2}y^2+(k-a-c-1+m)x^2y^{-2}.
\nonumber
\end{eqnarray} 
\end{enumerate}
\begin{proposition}
\label{Prop_SRGiii}
Let $X$ be a strongly regular graph with parameters
\begin{equation*}
(4\th^2-1,2\th^2,\th^2,\th^2),
\quad \text{for $\th \geq 2$.}
\end{equation*}
Let $W = I+xA(X)+yA(\overline{X})$ where
\begin{equation*}
x=-1 \quad \text{and} \quad y=\frac{2\th^2-3\pm \sqrt{4\th^2-5}\ i}{2(\th^2-1)}
\end{equation*}
or
\begin{equation*}
x=\frac{-2\th^2+1\pm \sqrt{4\th^2-1}\ i}{2\th^2}
\quad \text{and} \quad 
y=1.
\end{equation*}
Then $\nomw = \spn\{I, J\}$.  Consequently, $W$ is not of \Dita-type.
\end{proposition}
\proof
Suppose
\begin{equation*}
x=-1 \quad \text{and} \quad y=\frac{2\th^2-3\pm \sqrt{4\th^2-5}\ i}{2(\th^2-1)}.
\end{equation*}
Using Equation~(\ref{Eqn_a_bb}), the real part of $<\y{\a}{\b},\y{\g}{\a}>$ is
\begin{equation*}
\frac{-(4\th^2-5)(\th^2-3)}{2(\th^2-1)^2}.
\end{equation*}
Using
Equations~(\ref{Eqn_a_b}),  (\ref{Eqn_ba_b}) and (\ref{Eqn_ba_bb}),
the real part of $<\y{\a}{\b}, \y{\g}{\a}>$ is
\begin{equation*}
\frac{-(4\th^2-5)}{2(\th^2-1)}
\end{equation*}
for all the other vertex $\g$ in $V(X)\backslash\{\a,\b\}$.
If $\th \geq 2$ then $<\y{\a}{\b},\y{\g}{\a}> \neq 0$ for all
$\g \neq \a,\b$.  
By Lemma~\ref{Lem_Trivial}, $\nomw = \spn\{I, J\}$.

Suppose 
\begin{equation*}
x=\frac{-2\th^2+1\pm \sqrt{4\th^2-1}\ i}{2\th^2}
\quad \text{and} \quad 
y=1.
\end{equation*}
Using Equations~(\ref{Eqn_a_b}), (\ref{Eqn_a_bb}), (\ref{Eqn_ba_b}) and
(\ref{Eqn_ba_bb}), 
the real part of $<\y{\a}{\b}, \y{\g}{\a}>$ is
\begin{equation*}
\frac{-(4\th^2-1)}{2\th^2}
\end{equation*}
which is non-zero for $\g \neq \a,\b$ and for $\th \geq 1$.
By Lemma~\ref{Lem_Trivial}, $\nomw = \spn\{I, J\}$.
\qed

\begin{proposition}
\label{Prop_SRGv}
Let $X$ be a strongly regular graph with parameters
\begin{equation*}
(4\th^2+4\th+2,2\th^2+\th,\th^2-1,\th^2).
\end{equation*}
Let $W = I+xA(X)+yA(\overline{X})$ where
\begin{equation*}
x=\frac{-1\pm\sqrt{4\th^2(\th+1)^2-1}\ i}{2\th(\th+1)}
\quad \text{and} \quad 
y=x\inv=\frac{-1\mp\sqrt{4\th^2(\th+1)^2-1}\ i}{2\th(\th+1)}
\end{equation*}
Then $\nomw = \spn\{I, J\}$.  Consequently, $W$ is not of \Dita-type.
\end{proposition}
\proof
If $\g$ is adjacent to both $\a$ and $\b$ then,
by Equation~(\ref{Eqn_a_b}),
the real part of $<\y{\a}{\b}, \y{\g}{\a}>$ is
\begin{equation*}
\frac{-(\th^2+\th+1)(2\th^2+2\th+1)(2\th^2+2\th-1)}{2\th^3(1+\th)^3}.
\end{equation*}
If $\g$ is adjacent to $\a$ but not to $\b$ then
by Equation~(\ref{Eqn_a_bb}),
the real part of $<\y{\a}{\b}, \y{\g}{\a}>$ is
\begin{equation*}
\frac{-(2\th^2+2\th-1)(2\th^2+2\th+1)(\th^4+2\th^3-\th-1)}{2\th^4(1+\th)^4}.
\end{equation*}
If $\g$ is not adjacent to $\a$ then
by Equations~(\ref{Eqn_ba_b}) and (\ref{Eqn_ba_bb}),
the real part of $<\y{\a}{\b}, \y{\g}{\a}>$ is
\begin{equation*}
\frac{(2\th^2+2\th-1)(2\th^2+2\th+1)}{2\th^2(1+\th)^2}.
\end{equation*}
We have $<\y{\a}{\b},\y{\g}{\a}> \neq 0$, for $\g \neq \a, \b$ and for $\th \geq 1$.
By Lemma~\ref{Lem_Trivial}, $\nomw = \spn\{I, J\}$.
\qed

\section{Acknowledgement}
The author wishes to thank Chris Godsil for many useful discussions, and thank Akihiro Munemasa and Takuya Ikuta for telling her about the fourth complex Hadamard matrix in the 
adjacency algebra of the line graph of Petersen graph.


\end{document}